\documentclass[11pt]{amsart}

\usepackage[all]{xy}
\usepackage{amssymb}
\usepackage{amsmath}
\usepackage{amsmath}
\usepackage{mathrsfs}

\def\alg{\textsf{-alg}}
\def\-{\textrm{-}}

\def\B{\bigskip}
\def\BB{\bigskip \bigskip}

\def\r{\succ}
\def\l{\prec}
\def\g{\dashv }
\def\d{\vdash }

\def\m{\, \bot\,  }

\def\DD{\Delta}

\def\ss{\sigma}

\def\t{\otimes}

\def\Vtn{V^{\otimes n}}
\def\id{\textrm{id}}
\def\Sy{{\mathbb S}}
\def\Ai{A_{\infty}}

\def\tH{\underset{\textrm{H}}{\otimes}}

\renewcommand{\ne}{\nearrow}
\newcommand{\se}{\searrow}
\newcommand{\sw}{\swarrow}
\newcommand{\nw}{\nwarrow}
\newcommand{\east}{\succ}
\newcommand{\west}{\prec}
\newcommand{\north}{\wedge}
\newcommand{\south}{\vee}

\def\KK{{\mathbb{K}}}

\def\AAA{{\mathcal{A}}}
\def\CCC{{\mathcal{C}}}

\def\PP{{\mathcal{P}}}
\def\QQQ{{\mathcal{Q}}}
\def\RRR{{\mathcal{R}}}

\def\KKK{{\mathcal{K}}}

\def\gg{\mathfrak{g}}

\def\Id{\mathrm{Id }}

\def\sgn{\mathrm{sgn}}

\def\mono{\rightarrowtail}

\def\Vect{\mathop{\texttt{Vect}}}
\def\dgVect{\mathop{\texttt{dgVect}}}
\def\Ind{\mathop{\rm Ind }}

\def\dim{\mathop{\rm dim }}

\def\To{\overline {T}}

\def\ac{{}^{\scriptstyle \textrm{!`}}}
\def\pt{{\textrm{-}}}

\begin{document}

\author{G. W. Zinbiel}
\address{Zinbiel Institute of Mathematics (France)}
\email{gw.zinbiel@free.fr}
\title{Encyclopedia of types of algebras 2010}

\maketitle 
\begin{abstract} This is a cornucopia of types of algebras with some of their properties from the operadic point of view.
\end{abstract}
\vskip 1cm

\section*{Introduction} The following is a list of some types of algebras together with their properties under an operadic and homological point of view.
We keep the information to one page per type and we provide one reference as Ariadne's thread. More references are listed by the end of the paper.
We work over a fixed field $\KK$ though in many instances everything makes sense and holds over a commutative ground ring. The category of vector spaces over $\KK$ is denoted by $\Vect$. All tensor products are over $\KK$ unless otherwise stated.

The items of a standard page (which is to be found at the end of this introduction) are as follows.

\bigskip

Sometimes a given type appears under different names in the literature. The choice made in \texttt{Name} is, most of the time, the most common one (up to a few exceptions). The other possibilities appear under the item \texttt{Alternative}. 

The presentation given in \texttt{Definit.} is the most common one (with one exception). When others are used in the literature they are given in  \texttt{Alternative}. The item  \texttt{oper.} gives the generating operations. The item  \texttt{sym.} gives their symmetry properties, if any. The item  \texttt{rel.} gives the relation(s). They are supposed to hold for any value of the variables $x,y,z, \ldots $. If, in the presentation, only binary operations appear, then the type is said to be \emph{binary}. Analogously, there are \emph{ternary, $k$-ary, multi-ary} types.

If, in the presentation, the relations involve only the composition of two operations at a time (hence 3 variables in the binary case, 5 variables in the ternary case), then the type is said to be \emph{quadratic}.

For a given type of algebras $\PP$ the category of $\PP$-algebras is denoted by $\PP\alg$.
For each type there is defined a notion of \emph{free algebra}. By definition the free algebra of type $\PP$ over the vector space $V$ is an algebra denoted by $\PP(V)$ satisfying the following universal condition:

for any algebra $A$ of type $\PP$ and any linear map $\phi:V\to A$ there is a unique $\PP$-algebra morphism $\tilde \phi:\PP(V)\to A$ which lifts $\phi$. In other words the forgetful functor 
$$\PP\alg \longrightarrow \Vect$$
admits a left adjoint
$$\PP :\Vect \longrightarrow \PP\alg .$$
In all the cases mentioned here the relations involved in the presentation of the given type are (or can be made) multilinear. Hence the functor $\PP(V)$ is of the form (at least in characteristic zero),
$$\PP(V) = \bigoplus_{n\geq 1}\PP(n)\t _{\Sy_{n}} \Vtn\ ,$$
where $\PP(n)$ is some $\Sy_{n}$-module. The $\Sy_n$-module $\PP(n)$ is called the space of $n$-ary operations since for any algebra $A$ there is a map 
$$\PP(n)\t_{\Sy_{n}} A^{\t n} \to A.$$

The functor $\PP:\Vect \to \Vect$ inherits a monoid structure from the properties of the free algebra. Hence there exist transformations of functors $\iota: \Id \to \PP$ and $\gamma : \PP\circ \PP \to \PP$ such that $\gamma$ is associative and unital. The monoid $(\PP, \gamma , \iota)$ is called a \emph{symmetric operad}. 

The symmetric operad $\PP$ can also be described as a family of $\Sy_{n}$-modules $\PP(n)$ together with maps 
$$\gamma(i_{1},\ldots , i_{k}): \PP(k)\t\PP(i_{1})\t \cdots \t \PP( i_{k})\longrightarrow \PP(i_{1}+\cdots+ i_{k})$$
satisfying some  compatibility with the action of the symmetric group and satisfying the associativity property. 

If $\Sy_{n}$ is acting freely on $\PP(n)$, then $\PP(n)= \PP_{n}\t \KK[\Sy_{n}]$ where $\PP_{n}$ is some vector space, and $\KK[\Sy_{n}]$ is the regular representation. If, moreover,  the maps 
$\gamma(i_{1},\ldots , i_{n})$ are induced by maps 
$$\gamma_{i_{1},\ldots , i_{k}}: \PP_{k}\t\PP_{i_{1}}\t \cdots \t \PP_{ i_{k}}\longrightarrow \PP_{i_{1}+\cdots+ i_{k}} ,$$
then the operad $\PP$ comes a from a \emph{nonsymmetric operad} (abbreviated ns operad), still denoted $\PP$. 

For more terminology and details about algebraic operads we refer to \cite{JLLBV}.

 The generating series of the operad $\PP$  is defined as
 $$f^{\PP}(t):=\sum_{n\geq 1}\frac{\dim \PP(n)}{n!}\ t^n ,$$
 in the  binary case.  When dealing with a nonsymmetric  operad it becomes
 $$f^{\PP}(t):=\sum_{n\geq 1}\dim \PP_n\ t^n.$$

The Koszul duality theory of associative algebras has been extended to binary quadratic operads by Ginzburg and Kapranov, cf.\ \cite{GinzburgKapranov94}, then to quadratic operads by Fresse,  cf.\ \cite{Fresse06}. We give a conceptual treatment of this theory, together with applications in \cite{JLLBV}.  So, to any quadratic operad $\PP$, there is associated a quadratic \texttt{dual operad} denoted $\PP^!$. It is often a challenge to find a presentation of $\PP^!$ out of a presentation of $\PP$. One of the main results of the Koszul duality theory of operads is to show the existence of a natural differential map on the composite $\PP^{!*}\circ \PP$ given rise to the \emph{Koszul complex}. If it is acyclic, then $\PP$ is said to be \emph{Koszul}. One can show that, if $\PP$ is Koszul, then so is $\PP^!$. In this case the generating series
are inverse to each other for composition, up to sign, that is:
$$f^{\PP^!}(-f^{\PP}(-t))=t.$$

In the $k$-ary case we introduce the skew-generating series
 $$g^{\PP}(t):=\sum_{n\geq 1}(-1)^k\frac{\dim \PP((k-1)n+1)}{n!}\ t^{((k-1)n+1)} .$$
 If the operad $\PP$ is Koszul, then (cf.\ \cite{Vallette07})
 $$f^{\PP^!}(-g^{\PP}(-t))=t.$$

The items \texttt{Free alg.}, \texttt{rep.} $\PP(n)$ or $\PP_{n}$, $\dim \PP(n)$ or $\dim \PP_{n}$, and \texttt{Gen.~series} speak for themselves.

Koszulity of an operad implies the existence of a small chain complex to compute the (co)homology of a $\PP$-algebra.  When possible, the information on it is given in the item
 \texttt{Chain-cplx}. Moreover it permits us to construct the notion of $\PP$-algebra up to homotopy, whose associated operad, which is a differential graded operad, is denoted $\PP_{\infty}$.  
 
 The item  \texttt{Properties} lists the main features of the operad.  \emph{Set-theoretic} means that there is a set operad $\PP_{Set}$ (monoid in the category of $\Sy$-Sets) such that $\PP = \KK[\PP_{Set}]$. Usually this property can be read on the presentation of the operad: no algebraic sums. 
 
 In the item \texttt{Relationship} we list some of the ways to obtain this operad under some natural constructions like tensor product (Hadamard product) or Manin products (white {\LARGE $\circ$} or black {\LARGE $\bullet$}), denoted $\square$ and $\blacksquare$ in the nonsymmetric framework, cf.\cite{JLLBV}) for instance. We also list some of the most common functors to other types of algebras. Keep in mind that a functor $\PP \to \QQQ$ induces a functor $\QQQ\alg \to \PP\alg$ on the categories of algebras.
 
 Though we describe only algebras without unit, for some types there is a possibility of introducing an element 1 which is either a unit or a partial unit for some of the operations, see the discuission in \cite{Loday04b}. We indicate it in the item \texttt{Unit}.
 
 For binary operads the \emph{opposite type} consists in defining new operations by $x\cdot y = yx$, etc. If the new type is isomorphic to the former one, then the operad is said to be \emph{self-opposite}. When it is not the case, we mention whether the given type is called \emph{right} or \emph{left} in the item \texttt{Comment}.
 
In some cases the structure can be ``integrated". For instance Lie algebras are integrated into Lie groups (Lie third problem). If so, we indicate it in  the item \texttt{Comment}.

 In the item \texttt{Ref.} we indicate a reference where information on the operad and/or on the (co)homology theory can be obtained. It is not necessarily the first paper in which this type of algebras first appeared. For the three graces $As, Com, Lie$, the classical books Cartan-Eilenberg ``Homological Algebra'' and MacLane ``Homology'' are standard references. 

\B
Here is the list of the types included so far (with page number and letter K indicating that they are Koszul dual to each other):

\BB

\begin{tabular}{ l r | l r | l r | l r | l}
sample       & 6&
$As		$&  7& self-dual\\
$Com	$&  8&
$Lie	        $&  9& K\\
$Pois	$& 10& 
none          & 11& self-dual\\
$Leib	$& 12&
$Zinb	$& 13&K\\
$Dend	$& 14&
$Dias	$& 15&K\\
$PreLie	$& 16&
$Perm	$& 17&K\\
$Dipt	$& 18&
$Dipt^!     $& 19&K\\
$2as	        $& 20&
$2as^!     $& 21&K\\
$Tridend $& 22&
$Trias	$& 23&K\\
$PostLie	$& 24&
$ComTrias$& 25& K\\
$CTD	$& 26&
$CTD^!	$& 27&K\\
$Gerst      $& 28&
$BV          $& 29\\
$Mag	$& 30&
$Nil_{2}	$& 31&K\\
$ComMag$& 32&
$ComMag^{!}$& 33&K\\
$Quadri  $& 34&
$Quadri^!$& 35&K\\
$Dup       $& 36&
$Dup^!    $& 37&K\\
$As^{(2)}  $& 38&
$As^{\langle 2\rangle}$& 39\\
$L\-dend $& 40&
$Lie\-adm$& 41\\
$PreLiePerm$& 42&
$Altern	    $& 43\\
$Param1rel$& 44&
$MagFine  $& 45\\
$GenMag$& 46&
$NAP       $& 47\\
$Moufang$& 48&
$Malcev$& 49\\
$Novikov$& 50&
$Double Lie$& 51\\
$DiPreLie$& 52&
$Akivis$& 53\\
$Sabinin$& 54&
$Jordan\ triples$& 55\\
$ t\-As^{(3)}$& 56&
$p\-As^{(3)}$& 57\\
$LTS            $& 58&
$Lie\-Yamaguti$& 59\\
$Interchange$&60&
$HyperCom$& 61\\
$\Ai             $& 62&
$C_{\infty}$& 63\\
$L_{\infty}$& 64&
$Dend_{\infty}$&65\\
${\PP}_{\infty}$& 66&
$Brace     $& 67\\
$MB          $& 68&
$2Pois$& 69\\
$\Xi^{\pm}$& 70&
your own& 71\\
\end{tabular}
\vfill

\newpage

\noindent{\bf Notation}

We use the notation $\Sy_n$ for the symmetric group. Trees are very much in use in the description of operads. We use the following notation:

\medskip

-- $PBT_n$ is the set of planar binary rooted trees with $n-1$ internal vertices (and hence $n$ leaves). The number of elements of $PBT_n$ is the Catalan number 
$c_n=\frac{1}{n+1}\binom{2n}{n}$.

--$PT_n$  is the set of planar rooted trees with $n$ leaves, whose vertices have valency greater than $1+2$ (one root, at least 2 inputs). So we have $PBT_n\subset PT_n$. The number of elements of $PT_n$ is the super Catalan number, also called Schr\"oder number, denoted $C_n$.

\medskip

A planar binary rooted tree $t$ is completely determined by its right part $t^r$ and its left part $t^l$. More precisely $t$ is the grafting of $t^l$ and $t^r$:
$$ t = t^l\vee t^r\ .$$

\noindent{\bf Comments}

Many thanks to Walter Moreira for setting up a software which computes the first dimensions of the operad from its presentation.

We remind the reader that we can replace the symmetric monoidal category $\Vect$ by many other symmetric monoidal categories. So there are notions of graded algebras, differential graded algebras, twisted algebras, and so forth.
In the graded cases the Koszul sign rule is in order. Observe that there are also operads where the operations may have different degree (operad encoding Gerstenhaber algebras for instance).
\bigskip

We end this paper with a tableau of integer sequences appearing in this document.

\bigskip

This list of types of algebras is not as encyclopedic as the title suggests. We put only the types which are defined by a finite number of generating operations and whose relations are multilinear. Moreover we only put those which have been used some way or another. You will not find the ``restricted types'' (like divided power algebras), nor bialgebras.
\B

Please report any error or comment to: \qquad \texttt {gw.zinbiel@free.fr} 
\bigskip

We plan to update this encyclopedia every now and then. 

\newpage

\section{Type of algebras }\index{algebra}

\begin{tabular}{ll}
\texttt{Name}                   & {\bf   Most common terminology}\\
&\\
\texttt{Notation}               & my favorite notation for the operad (generic notation: $\PP$)\\
&\\
\texttt{Def.\ oper.}         & list of the generating operations\\
\texttt{\hskip 1.2cm sym.}  & their symmetry if any\\
\texttt{\hskip 1.3cm rel.}    & the relation(s)\\
&\\
\texttt{Free alg.}& the free algebra as a functor in $V$\\
&\\
\texttt{rep.~$\PP(n)$}      & the $\Sy_n$-representation $\PP(n)$ and/or the space $\PP_n$ if nonsymmetric\\
&\\
\texttt{$\dim \PP(n)$}       & the series (if  close formula  available), the list of the 7 first numbers\\
&beginning at $n=1$\\
&\\
\texttt{Gen.~series}& close formula for  $f^{\PP}(t)=\sum_{n\geq 1}\frac{\dim \PP(n)}{n!}t^n$ when available\\
&\\
\texttt{Dual operad}& the Koszul dual operad \\
&\\
\texttt{Chain-cplx}&Explicitation of the chain complex, if not too complicated\\
&\\
\texttt{Properties}          &  among: nonsymmetric, binary, quadratic, set-theoretic, ternary, multi-ary,\\
& cubic, Koszul.\\
&\\
\texttt{Alternative}         & alternative terminology, and/or notation, and/or presentation.\\
& \\
\texttt{Relationship}      &some of the relationships with other operads, either under some\\
& construction like symmetrizing, Hadamard product, Manin products, or\\
& under the existence of functors\\
&\\
\texttt{Unit}     &whether one can assume the existence of a unit (or partial unit)\\
&\\
\texttt{Comment}          & whatever needs to be said which does not fit into the other items\\
&\\
\texttt{Ref.}&a reference, usually dealing with the homology of the $\PP$-algebras\\
&(not necessarily containing all the results of this page)\\
\end{tabular}

\newpage


\begin{tabular}{ll}
\texttt{Name}& {\bf Associative algebra}\index{algebra!associative}\\
&\\
\texttt{Notation}& $As$ (as nonsymmetric operad)\\
&$Ass$ (as symmetric operad)\\
&\\
\texttt{Def.\ oper.}& $xy$, operadically $\mu$\\
\texttt{\hskip 1.2cm sym.}  &\\
\texttt{\hskip 1.3cm rel.} &$(xy)z=x(yz)$, operadically $\mu\circ_{1}\mu = \mu\circ_{2}\mu$ (associativity)\\
&\\
\texttt{Free alg.}& $As(V)= \overline{T}(V)=\oplus_{n\ge 1}V^{\t n}$\quad  tensor algebra\\
& (noncommutative polynomials)\\
&$(x_1\ldots x_p)(x_{p+1}\ldots x_{p+q})= x_1 \ldots x_{p+q}$ (concatenation)\\
&\\
\texttt{rep.~$\PP(n)$}& $As(n)=\KK[\Sy_n]$ (regular representation),\quad $As_n = \KK$ \\
&\\
\texttt{$\dim \PP(n)$}
&$1, 2, 6, 24, 120, 720, 5040, \ldots, n!, \ldots $\\
&\\
\texttt{Gen.~series}& $f^{As}(t)=\frac{t}{1-t}$\\
&\\
\texttt{Dual operad}& $As^!=As$\\
&\\
\texttt{Chain-cplx}& non-unital Hochschild complex, $C_n^{As}(A)= A^{\t n}$ \\
&$b'(a_1, \ldots , a_n) = \sum_{i=1}^{i=n-1}(-1)^{i-1}(a_1,\ldots, a_{i}a_{i+1}, \ldots , a_{n})$\\
&important variation: cyclic homology\\
&\\
\texttt{Properties}& binary, quadratic, ns, set-theoretic, Koszul, self-dual.\\
&\\
\texttt{Alternative}& Associative algebra often simply called \emph{algebra}.\\
&Can be presented with commutative operation $x\cdot y:= xy+yx$\\
& and anti-symmetric operation $[x,y]=xy-yx$ satisfying\\
& {$\begin{cases}
       [x \cdot y, z] = x\cdot[y,z] + [x,z]\cdot y, \\
       (x\cdot y)\cdot z -x\cdot (y \cdot z) = [x,[y,z]]\ .\\
      \end{cases}$
  }\\
    
& (Livernet and Loday, unpublished)\\
&\\
\texttt{Relationship}& $Ass\alg \to Lie\alg,\  [x,y]= xy-yx$,\\
& $Com\alg \to Ass\alg$ (inclusion), and many others\\
&\\
\texttt{Unit}& $1x=x=x1$\\
&\\
\texttt{Comment}& one the three graces\\
\end{tabular}

\newpage


\begin{tabular}{ll}
\texttt{Name}& {\bf   Commutative algebra}\index{algebra!commutative}\\
&\\
\texttt{Notation}& $Com$\\
&\\
\texttt{Def.\ oper.}& $xy$\\
\texttt{\hskip 1.2cm sym.}  & $xy=yx $ (commutativity)\\
\texttt{\hskip 1.3cm rel.} &$(xy)z=x(yz)$ (associativity)\\
&\\
\texttt{Free alg.}&$Com(V)=\overline{S}(V)$ (polynomials)   \\
&if $V=\KK x_{1}\oplus \cdots \oplus \KK x_{n}$, then $\KK 1\oplus Com(V) = \KK [x_{1}, \ldots , x_{n}]$\\
\texttt{rep.~$\PP(n)$}& $Com(n)= \KK$ (trivial representation)\\
&\\
\texttt{$\dim \PP(n)$}&$1, 1, 1, 1, 1, 1, 1, \ldots, 1, \ldots $ \\
&\\
\texttt{Gen.~series}& $f^{Com}(t)=\exp(t) - 1$\\
&\\
\texttt{Dual operad}& $Com^!=Lie$\\
&\\
\texttt{Chain-cplx}& Harrison complex in char. 0, Andr\'e-Quillen cplx in general  \\
&\\
\texttt{Properties}&  binary, quadratic, set-theoretic, Koszul.\\
&\\
\texttt{Alternative}&sometimes called associative and commutative algebra\\
& other notation $Comm$\\
&\\
\texttt{Relationship}&$Com\alg \to As\alg$, \\
&$Zinb\alg \to Com\alg$\\
&\\
\texttt{Unit}& $1x=x=x1$\\
&\\
\texttt{Comment}& one of the three graces\\

\end{tabular}

\newpage


\begin{tabular}{ll}
\texttt{Name}& {\bf   Lie algebra}\index{algebra!Lie}\\
&\\
\texttt{Notation}& $Lie$\\
&\\
\texttt{Def.\ oper.}& $[x,y]$ (bracket)\\
\texttt{\hskip 1.2cm sym.}  & $[x,y]=-[y,x] $ (anti-symmetry)\\
\texttt{\hskip 1.3cm rel.} &$[[x,y],z]= [x,[y,z]] + [[x,z],y]$ (Leibniz relation)\\
&\\
\texttt{Free alg.}&$Lie(V)=$ subspace of the tensor algebra $T(V)$\\
& generated by $V$ under the bracket \\
&\\
\texttt{rep.~$\PP(n)$}& $Lie(n)=\Ind_{C_n}^{\Sy_n}(\sqrt[n]{1})$\\
&\\
\texttt{$\dim \PP(n)$}
&$1, 1, 2, 6, 24, 120, 720,  \ldots , (n-1)!, \ldots$\\
&\\
\texttt{Gen.~series}& $f^{Lie}(t)=-\log(1-t)$\\
&\\
\texttt{Dual operad}& $Lie^!=Com$\\
&\\
\texttt{Chain-cplx}& Chevalley-Eilenberg complex  $C_n^{Lie}(\gg)= \Lambda^n\gg$   \\
&$d(x_1\wedge \ldots \wedge x_n)=$\\
&\qquad $ \sum_{1\leq i<j\leq n}(-1)^j(x_1\wedge  \ldots\wedge  [x_i,x_j]\wedge \ldots \wedge \widehat{x_j}\wedge  \ldots \wedge  x_n)$\\
&\\
\texttt{Properties}&  binary, quadratic, Koszul.\\
&\\
\texttt{Alternative}&The relation is more commonly written as the Jacobi identity:\\
&$[x,[y,z]] + [y,[z,x]]+ [z,[x,y]] =0$\\
&\\
\texttt{Relationship}&$As\alg \to Lie\alg$, \\
&$Lie\alg \to Leib\alg$, \\
&$PreLie\alg \to Lie\alg$\\
&\\
\texttt{Unit}& no\\
&\\
\texttt{Comment}& one of the three graces. Named after Sophus Lie.\\
& Integration: Lie groups.\\

\end{tabular}

\newpage


\begin{tabular}{ll}
\texttt{Name}& {\bf   Poisson algebra}\index{algebra!Poisson}\\
&\\
\texttt{Notation}& $Pois$\\
&\\
\texttt{Def.\ oper.}& $xy, \{x,y\}$\\
\texttt{\hskip 1.2cm sym.}  & $xy = yx, \{x,y\}=-\{y,x\} $\\
\texttt{\hskip 1.3cm rel.} & {$\begin{cases}
       \{\{x,y\},z\}= \{x,\{y,z\}\} + \{\{x,z\},y\},\\
       \{x y, z\} = x\{y,z\} + \{x,z\}y, \\
       (x y) z -x (y  z) = 0\ .\\
      \end{cases}$
  }\\
&\\
\texttt{Free alg.}&$Pois(V)\cong \overline{T}(V)$ (tensor module, iso as Schur functors)\\
&\\
\texttt{rep.~$\PP(n)$}& $Pois(n)\cong\KK[\Sy_n]$ (regular representation)\\
&\\
\texttt{$\dim \PP(n)$}
&$1, 2!, 3!, 4!, 5!, 6!, 7!, \ldots , n!, \ldots $ \\
&\\
\texttt{Gen.~series}& $f^{Pois}(t)=\frac{t}{1-t}$\\
&\\
\texttt{Dual operad}& $Pois^!=Pois$\\
&\\
\texttt{Chain-cplx}& Isomorphic to the total complex of a certain bicomplex cons-\\
&tructed from the action of the Eulerian idempotents  \\
&\\
\texttt{Properties}&  binary, quadratic, quasi-regular, Koszul, self-dual.\\
&\\
\texttt{Alternative}&Can be presented with one operation $x*y$ satisfying the relation\\
&$(x*y)*z =x*(y*z)+$\\
&\qquad $   \frac{1}{3}\big(+x*(z*y) - z*(x*y) -y*(x*z)+y*(z*x)\big)$\\
&\\
\texttt{Relationship}&$Pois\alg \rightleftarrows Lie\alg$,\quad  $Pois\alg \rightleftarrows Com\alg$,\\
&\\
\texttt{Unit}& $1x=x=x1, [1,x]=0=[x,1]$\\
&\\
\texttt{Comment}&Named after Sim\'eon Poisson.  \\
&\\
\texttt{Ref.}
&\cite{Fresse06} B.\ Fresse, \emph{Th\' eorie des op\' erades de Koszul et homologie des}\\
& \emph{ alg\`ebres de Poisson,} Ann.Math. Blaise Pascal 13 (2006), 237--312. \\
\end{tabular}

\newpage

This page is inserted so that, in the following part, an operad and its dual appear on page $2n$ and $2n+1$ respectively. Since $Pois$ is self-dual there is no point to write a page about its dual.

\vskip2cm
Let us take the opportunity to mention that if $A$ is a $\PP$-algebra and $B$ is a $\PP^!$-algebra, then the tensor product $A \t B$ inherits naturally a structure of Lie algebra. If $\PP$ is nonsymmetric, then so is $\PP^!$, and  $A \t B$ is in fact an associative algebra.

In some  cases (like the Leibniz case for instance),  $A \t B$ is a pre-Lie algebra.
\newpage


\begin{tabular}{ll}
\texttt{Name}& {\bf   Leibniz algebra}\index{algebra!Leibniz}\\
&\\
\texttt{Notation}& $Leib$\\
&\\
\texttt{Def.\ oper.}         & $[x,y]$\\
\texttt{\hskip 1.2cm sym.}  & $ $\\
\texttt{\hskip 1.3cm rel.}    & $[[x,y],z]= [x,[y,z]] + [[x,z],y]$ (Leibniz relation)\\
&\\
\texttt{Free alg.}& $Leib(V)\cong \To(V)$ (reduced tensor module, iso as Schur functors)\\
&\\
\texttt{rep.~$\PP(n)$}      & $Leib(n)=\KK[\Sy_n]$ (regular representation)\\
&\\
\texttt{$\dim \PP(n)$}   
&1, 2!, 3!, 4!, 5!, 6!, 7!, \ldots , n!, \ldots \\
&\\
\texttt{Gen.~series}& $f^{Leib}(t)=\frac{t}{1-t}$\\
&\\
\texttt{Dual operad}& $Leib^!=Zinb$\\
&\\
\texttt{Chain-cplx}& $C_n^{Leib}(\gg)= \gg^{\t n}$   \\
&$d(x_1,\ldots , x_n)= \sum_{1\leq i<j\leq n}(-1)^j(x_1, \ldots, [x_i,x_j],\ldots ,\widehat{x_j}, \ldots , x_n)$\\
&\\
\texttt{Properties}          &  binary, quadratic, quasi-regular, Koszul.\\
&\\
\texttt{Alternative}         &   Sometimes improperly called Loday algebra.  \\
& \\
\texttt{Relationship}      &$Leib = Perm$ {\LARGE $\circ$} $Lie$ (Manin white product), see \cite{Vallette08}\\
&$Lie\alg \to Leib\alg ,\ Dias\alg \to Leib\alg , \ Dend\alg \to Leib\alg $\\
&\\
\texttt{Unit}     & no\\
&\\
\texttt{Comment}          & Named after G.W.~Leibniz. This is the \emph{left} Leibniz algebra. \\
&The opposite type is called \emph{right} Leibniz algebra.\\
&Integration : ``coquecigrues'' ! see for instance \cite{Covez10}\\
&\\
\texttt{Ref.} &\cite{Loday93} J.-L. Loday, \emph{Une version non commutative des alg\`ebres}\\
& \emph{ de Lie: les alg\`ebres de Leibniz.}\\
& Enseign.~Math.~(2) 39 (1993), no. 3-4, 269--293.\\

\end{tabular}

\newpage


\begin{tabular}{ll}
\texttt{Name}                       & {\bf  Zinbiel algebra}\index{algebra!Zinbiel}\\
&\\
\texttt{Notation}                   & $Zinb$\\
&\\
\texttt{Def.\ oper.}            & $x\cdot y$\\
\texttt{\hskip 1.2cm sym.}     & $ $\\
\texttt{\hskip 1.3cm rel.}       & $(x\cdot y)\cdot z = x\cdot (y\cdot z) + x\cdot (z\cdot y)\ $ (Zinbiel relation)\\
&\\
\texttt{Free alg.}  &  $Zinb(V)=\overline{T}(V),\quad  \cdot=\mathrm{halfshuffle}$ \\
& $x_1\ldots x_p \cdot x_{p+1}\ldots x_{p+q} = x_1 \mathrm{sh}_{p-1, q}(x_2\ldots x_p ,x_{p+1}\ldots x_{p+q} )$\\
&\\
\texttt{rep.~$\PP(n)$}      & $Zinb(n)\cong \KK[\Sy_n]$ (regular representation)\\
&\\
\texttt{$\dim \PP(n)$}   
&1, 2!, 3!, 4!, 5!, 6!, 7!, \ldots , n!, \ldots \\
&\\
\texttt{Gen.~series}& $f^{Zinb}(t)=\frac{t}{1-t}$\\
&\\
\texttt{Dual operad}& $Zinb^!=Leib$\\
&\\
\texttt{Chain-cplx}&   \\
&\\
\texttt{Properties}            &  binary, quadratic, quasi-regular, Koszul.\\
&\\
\texttt{Alternative}           & $Zinb= ComDend$ (commutative dendriform algebra)  \\
&Previously called \emph{dual Leibniz algebra}. \\
& \\
\texttt{Relationship}         & $Zinb = PreLie \bullet Com$, see \cite{Vallette08}, $Zinb=ComDend$\\
&$Zinb\alg \to Com\alg,\quad  xy=x\cdot y + y \cdot x$,\\
&  $Zinb\alg \to Dend\alg,\quad  x\l y=x\cdot y= y \r x$\\
&\\
\texttt{Unit}        & $1\cdot x=0,\quad  x\cdot 1= x$\\
&\\
\texttt{Comment}             & symmetrization of the dot product gives, not only\\ 
&a commutative alg., but in fact a \emph{divided power algebra}.\\
&Named after G.W.\ Zinbiel. This is right Zinbiel algebra.\\
&\\
\texttt{Ref.} &\cite{Loday95}  J.-L. Loday,  \emph{Cup-product for Leibniz cohomology and}\\
&\emph{dual Leibniz algebras.} Math.\ Scand.\ 77 (1995), no. 2, 189--196.\\
\end{tabular}

\newpage


\begin{tabular}{ll}
\texttt{Name}                   & {\bf Dendriform  algebra}\index{algebra!dendriform}\\
&\\
\texttt{Notation}               & $Dend$\\
&\\
\texttt{Def.\ oper.}         & $x\l y,\  x\r y$  (left and right operation)\\
\texttt{\hskip 1.2cm sym.}  & $ $\\
\texttt{\hskip 1.3cm rel.}    & 
$\left\{\begin{array}{rcl}
(x\l y)\l z &=& x\l (y\l z)+ x\l (y\r z)   , \\
(x\r y)\l z &=& x\r (y\l z) , \\
(x\l y)\r z + (x\r y)\r z&=& x\r (y\r z).
\end{array}
\right .
$
\\
&\\
\end{tabular}

\newpage


\begin{tabular}{ll}
\texttt{Name}                   & {\bf Diassociative  algebra}\index{algebra!diassociative}\\
&\\
\texttt{Notation}               & $Dias$\\
&\\
\texttt{Def.\ oper.}         & $x\g y, \ x\d y$ (left and right operation)\\
\texttt{\hskip 1.2cm sym.}  & $ $\\
\texttt{\hskip 1.3cm rel.}    & $\left\{\begin{array}{rcl}
(x\g y)\g z &=& x\g (y\g z), \\
(x\g y)\g z &=& x\g (y\d z), \\
(x\d y)\g z &=& x\d (y\g z) , \\
(x\g y)\d z &=& x\d (y\d z),\\
(x\d y)\d z &=& x\d (y\d z).
\end{array}
\right .
$
\\
&\\
\texttt{Free alg.}&$Dias(V)=\bigoplus_{n\geq 1}(\underbrace{\Vtn\oplus \cdots \oplus \Vtn}_{n\ \textrm{ copies}})$ \\
&noncommutative polynomials with one variable marked  \\
&\\
\texttt{rep.~$\PP(n)$}      & $Dias(n)=\KK^n\t \KK[\Sy_n], \quad Dias_n= \KK^n$\\
&\\
\texttt{$\dim \PP(n)$}       &$1, 2\times 2!, 3\times 3!, 4\times 4!, 5\times 5!, 6\times 6!, 7\times 7!, \ldots , n\times n!, \ldots,$\\
&\\
\texttt{Gen.~series}& $ f^{Dias}(t)= \frac{t}{(1-t)^2}$\\
&\\
\texttt{Dual operad}& $Dias^!=Dend$\\
&\\
\texttt{Chain-cplx}& see ref.  \\
&\\
\texttt{Properties}          &  binary, quadratic, ns, set-theoretic, Koszul.\\
&\\
\texttt{Alternative}         & Also called associative dialgebras, or for short, dialgebras.  \\
& \\
\texttt{Relationship}      &$Dias = Perm$ {\LARGE $\circ$} $As = Perm \tH As$\\
&$As\alg \to Dias\alg,\quad  x\g y := xy =: x\d y$\\
&$Dias\alg \to Leib\alg,\quad  [x,y]:= x\g y - x \d y$ \\
&\\
\texttt{Unit}   &Bar-unit: $x\g 1=x=1\d x,\ 1\g x= 0 = x\d 1$\\
&\\
\texttt{Comment}          & \\
&\\
\texttt{Ref.}&\cite{Loday01} J.-L.\ Loday, \emph{Dialgebras,} \\
&Springer Lecture Notes in Math. 1763 (2001), 7-66.\\
\end{tabular}

\newpage


\begin{tabular}{ll}
\texttt{Name}                   & {\bf Pre-Lie  algebra}\index{algebra!pre-Lie}\\
&\\
\texttt{Notation}               & $PreLie$\\
&\\
\texttt{Def.\ oper.}         & $x\circ y$\\
\texttt{\hskip 1.2cm sym.}  & $ $\\
\texttt{\hskip 1.3cm rel.}    & $(x\circ y) \circ z - x\circ (y\circ z) = (x\circ z) \circ y - x\circ (z\circ y) $\\
&\\
\texttt{Free alg.}& $PreLie(V)=\{\textrm{rooted trees labeled by elements of $V$}\}$ \\
&\\
\texttt{rep.~$\PP(n)$}      & $PreLie(n)=\KK[\{\textrm{rooted trees, vertices labeled by } 1, \ldots ,n\}]$ \\
&\\
\texttt{$\dim \PP(n)$}       & $1, 2, 9, 64, 625, 1296, 117649, \ldots, n^{n-1}, \ldots $\\
&\\
\texttt{Gen.~series}& $ f^{PreLie}(t)= y$ which satisfies $y=t\ \exp(y)$\\
&\\
\texttt{Dual operad}& $PreLie^!=Perm$\\
&\\
\texttt{Chain-cplx}&see ref.   \\
&\\
\texttt{Properties}          &  binary, quadratic, Koszul.\\
&\\
\texttt{Alternative}         &The relation is $as(x,y,z) = as(x,z,y)$.\\
& \\
\texttt{Relationship}      &$PreLie\alg \to Lie\alg,\quad [x,y] = x\circ y - y \circ x$ \\
& $Dend\alg \to PreLie\alg,\quad x\circ y := x\l y - y\r x$\\
& $Brace\alg \to PreLie\alg,$\quad forgetful functor\\
&\\
\texttt{Unit}     &$1\circ x = x = x\circ 1$\\
&\\
\texttt{Comment}          &This is \emph{right} pre-Lie algebra, also called \emph{right-symmetric algebra},\\
& or \emph{Vinberg algebra}\index{algebra!Vinberg}. The opposite type is left symmetric.\\
&First appeared in \cite{Gerstenhaber63, Vinberg63}.\\
&\\
\texttt{Ref.}&\cite{ChapotonLivernet01}  F.\ Chapoton, M.\ Livernet, \emph{Pre-Lie algebras and the rooted}\\
&\emph{trees operad.} Internat. Math. Res. Notices 2001, no. 8, 395--408.\\
\end{tabular}

\newpage


\begin{tabular}{ll}
\texttt{Name}                   & {\bf Perm  algebra}\index{algebra!permutative}\\
&\\
\texttt{Notation}               & $Perm$\\
&\\
\texttt{Def.\ oper.}         & $xy$\\
\texttt{\hskip 1.2cm sym.}  & $$\\
\texttt{\hskip 1.3cm rel.}    & $(xy)z=x(yz)=x(zy)$\\
&\\
\texttt{Free alg.}& $Perm(V)=V\t S(V)$\\
&\\
\texttt{rep.~$\PP(n)$}      & $Perm(n)=\KK^n$\\
&\\
\texttt{$\dim \PP(n)$}   
&$1, 2, 3, 4, 5, 6, 7,\ldots, n, \ldots$\\
&\\
\texttt{Gen.~series}& $ f^{Perm}(t)= t\exp(t)$\\
&\\
\texttt{Dual operad}& $Perm^!=PreLie $\\
&\\
\texttt{Chain-cplx}&   \\
&\\
\texttt{Properties}          &  binary, quadratic, set-theoretic, Koszul.\\
&\\
\texttt{Alternative}         &   \\
& \\
\texttt{Relationship}      &$Perm=ComDias$\\
& $Com\alg \to Perm\alg$,\\
&  $NAP\alg \to Perm\alg$,\\
&$Perm\alg \to Dias\alg$\\
&\\
\texttt{Unit}     &no, unless it is a commutative algebra\\
&\\
\texttt{Comment}          & \\
&\\
\texttt{Ref.}&\cite{Chapoton01a} F.\ Chapoton, \emph{Un endofoncteur de la cat\' egorie des} \\
& \emph{op\'erades.} Dialgebras \& related operads, 105--110,\\
& Lect. Notes in Math., 1763, Springer, 2001.\\
\end{tabular}

\newpage


\begin{tabular}{ll}
\texttt{Name}                   & {\bf Dipterous  algebra}\index{algebra!dipterous}\\
&\\
\texttt{Notation}               & $Dipt$\\
&\\
\texttt{Def.\ oper.}         & $x* y, \ x\l y $\\
\texttt{\hskip 1.2cm sym.}  & $ $\\
\texttt{\hskip 1.3cm rel.}    & 
$\left\{\begin{array}{rcl}
(x*y)* z &=& x* (y*z)\quad  \textrm{(associativity)}  \\
(x\l y)\l z &=& x\l (y* z)\quad  \textrm{(dipterous relation)}
\end{array}
\right .
$
\\
&\\
\texttt{Free alg.}& $Dipt(V)=\bigoplus_{n\geq 1}( \KK[PT_n]\oplus \KK[PT_n])\t V^{\t n}$, \\
& for $n=1$ the two copies of $PT_1$ are identified\\
&\\
\texttt{rep.~$\PP(n)$}      & $Dipt(n)= ( \KK[PT_n]\oplus \KK[PT_n])\t \KK[\Sy_n], n\geq 2$,\\ 
&$Dipt_n=  \KK[PT_n]\oplus \KK[PT_n], n\geq 2$\\
&\\
\texttt{$\dim \PP(n)$}   &$1, 2\times 2!, 6\times 3!, 22\times 4!, 90\times 5!, 394\times 6!, 1806\times 7!, \ldots , 2C_{n}\times n!, \ldots$\\
 &where $C_n$ is the Schr\"oder number: $\sum_{n\geq 1}C_n t^n = \frac{1+t-\sqrt{1-6t+t^2}}{4}$\\

&\\
\texttt{Gen.~series}& $f^{Dipt}(t)=  \frac{1-t-\sqrt{1-6t+t^2}}{2}$\\
&\\
\texttt{Dual operad}& $Dipt^! $\\
&\\
\texttt{Chain-cplx}& Isomorphic to the total complex of a certain explicit bicomplex  \\
&\\
\texttt{Properties}          &  binary, quadratic, ns, Koszul.\\
&\\
\texttt{Alternative}         &``dipt\`ere" in French  \\
& \\
\texttt{Relationship} 
&$Dend\alg \to Dipt\alg,\quad  x*y:=x\l y + y \r x$,\\
&a variation: replace the dipterous relation by  \\
&$(x\l y)*z + (x*y)\l z = x\l (y* z)+ x*(y\l z)$ \\
& to get $Hoch$-algebras, see \cite{Leroux11}. Same properties.\\
&\\
\texttt{Unit}       & $1\l x = 0 ,\ x\l 1 = x ,\quad
1*x = x=x* 1$.\\
&\\
\texttt{Comment}          &dipterous = 2-fold  in greek (the free algebra has 2  planar trees copies) \\
&\\
\texttt{Ref.} &\cite{LodayRonco03} J.-L.~Loday, M. Ronco, \emph{Alg\`ebres de Hopf colibres},\\
& C.R.Acad. Sci. Paris t. 337, Ser. I (2003), 153 -158.\\
\end{tabular}

\newpage


\begin{tabular}{ll}
\texttt{Name}                   & {\bf Dual dipterous  algebra}\index{algebra!dual dipterous}\\
&\\
\texttt{Notation}               & $Dipt^! $\\
&\\
\texttt{Def.\ oper.}         & $x\g y, \ x* y$\\
\texttt{\hskip 1.2cm sym.}  & $ $\\
\texttt{\hskip 1.3cm rel.}    & $\left\{\begin{array}{rcl}
(x* y)* z &=& x* (y*z), \\
(x\g y)\g z &=& x\g (y* z), \\
(x*y)\g z &=&0 , \\
(x\g y)* z &=& 0,\\
0&=& x* (y\g z)\\
0 &=& x\g (y\g z).
\end{array}
\right .
$
\\
&\\
\texttt{Free alg.}&$Dipt^! (V)=\overline{T}(V)\oplus \overline{T}(V)$\\
& \\
&\\
\texttt{rep.~$\PP(n)$}      & $Dipt^! (n)=\KK^2\t \KK[\Sy_n], \ Dipt^!_n = \KK^2, n\geq 2 $\\
&\\
\texttt{$\dim \PP(n)$}       &$1, 2\times 2!, 2\times 3!, 2\times 4!, 2\times 5!, 2\times 6!, 2\times 7!, \ldots , 2\times n!, \ldots,$\\
&\\
\texttt{Gen.~series}& $ f^{Dipt^!}(t)= \frac{t+t^2}{1-t}$\\
&\\
\texttt{Dual operad}& $Dipt^{!!}=Dipt$\\
&\\
\texttt{Chain-cplx}& see ref.  \\
&\\
\texttt{Properties}          &  binary, quadratic, ns, Koszul.\\
&\\
\texttt{Alternative}         &   \\
& \\
\texttt{Relationship}      
&$As\alg \to Dipt^! \alg,\quad  x\g y := xy =: x\d y$\\
&\\
\texttt{Unit}   &No\\
&\\
\texttt{Comment}          & \\
&\\
\texttt{Ref.} &\cite{LodayRonco03} J.-L.~Loday, M. Ronco, \emph{Alg\`ebres de Hopf colibres}, \\
&C.R.Acad. Sci. Paris t. 337, Ser. I (2003), 153 -158.\\
\end{tabular}

\newpage


\begin{tabular}{ll}
\texttt{Name}                   & {\bf Two-associative  algebra}\index{algebra!Two-associative}\\
&\\
\texttt{Notation}               & $2as$\\
&\\
\texttt{Def.\ oper.}         & $x* y, \ x\cdot y, $\\
\texttt{\hskip 1.2cm sym.}  & $ $\\
\texttt{\hskip 1.3cm rel.}    & 
$\left\{\begin{array}{rcl}
(x*y)* z &=& x* (y*z) , \\
(x\cdot y)\cdot  z &=& x\cdot (y\cdot  z) .
\end{array}
\right .
$
\\
&\\
\texttt{Free alg.}& $2as(V)=\bigoplus_{n\geq 1}( \KK[T_n]\oplus \KK[T_n])\t V^{\t n}$, where $T_n=$ planar trees \\
& for $n=1$ the two copies of $T_1$ are identified\\
&\\
\texttt{rep.~$\PP(n)$}      & $2as(n)= ( \KK[T_n]\oplus \KK[T_n])\t \KK[\Sy_n] $, so $Dipt_n=  \KK[T_n]\oplus \KK[T_n] , n\geq 2$\\
&\\
\texttt{$\dim \PP(n)$}   &$1, 2\times 2!, 6\times 3!, 22\times 4!, 90\times 5!, 394\times 6!, 1806\times 7!, \ldots , 2C_{n}\times n!, \ldots$\\
 &where $C_n$ is the Schr\"oder number: $\sum_{n\geq 1}C_n t^n = \frac{1+t-\sqrt{1-6t+t^2}}{4}$\\

&\\
\texttt{Gen.~series}& $f^{2as}(t)=  \frac{1-t-\sqrt{1-6t+t^2}}{2}$\\
&\\
\texttt{Dual operad}& $2as^! $\\
&\\
\texttt{Chain-cplx}& Isomorphic to the total complex of a certain explicit bicomplex  \\
&\\
\texttt{Properties}          &  binary, quadratic, ns, set-theoretic,  Koszul.\\
&\\
\texttt{Alternative}         &  \\
& \\
\texttt{Relationship}  & $2as\alg \to Dup\alg,\quad  2as\alg \to B_{\infty}\alg$ \\
& \\
&\\
\texttt{Unit}       & $1\cdot x =  x= x\cdot 1 ,\quad
1*x = x=x* 1$.\\
&\\
\texttt{Comment}          & \\
&\\
\texttt{Ref.} &\cite{LodayRonco06} J.-L.~Loday, M. Ronco, \emph{On the structure of cofree Hopf}\\
& \emph{algebras}, J. reine angew. Math. 592 (2006), 123--155.\\
\end{tabular}

\newpage


\begin{tabular}{ll}
\texttt{Name}                   & {\bf Dual 2{\textrm -}associative  algebra}\index{algebra!dual two associative}\\
&\\
\texttt{Notation}               & $2as^! $\\
&\\
\texttt{Def.\ oper.}         & $x\cdot  y, \ x* y$\\
\texttt{\hskip 1.2cm sym.}  & $ $\\
\texttt{\hskip 1.3cm rel.}    & $\left\{\begin{array}{rcl}
(x* y)* z &=& x* (y*z), \\
(x\cdot  y)\cdot  z &=& x\cdot  (y\cdot z), \\
(x\cdot y)* z &=&0 , \\
(x* y)\cdot  z &=& 0,\\
0&=& x* (y\cdot  z)\\
0 &=& x\cdot  (y* z).
\end{array}
\right .
$
\\
&\\
\texttt{Free alg.}&$2as^! (V)=V\oplus \oplus_{n\geq 2}\ \Vtn \oplus \Vtn $ \\
& \\
&\\
\texttt{rep.~$\PP(n)$}      & $Dipt^! (n)= \KK[\Sy_{n}]\oplus  \KK[\Sy_{n}],\quad Dipt^!_{n}= \KK\oplus  \KK, n\geq 2.$\\
&\\
\texttt{$\dim \PP(n)$}       &$1, 2\times 2!, 2\times 3!, 2\times 4!, 2\times 5!, 2\times 6!, 2\times 7!, \ldots , 2\times n!, \ldots $\\
&\\
\texttt{Gen.~series}& $ f^{2as^!}(t)= \frac{t+t^2}{1-t}$\\
&\\
\texttt{Dual operad}& $(2as^!)^!=2as$\\
&\\
\texttt{Chain-cplx}& see ref.  \\
&\\
\texttt{Properties}          &  binary, quadratic, Koszul.\\
&\\
\texttt{Alternative}         &   \\
& \\
\texttt{Relationship}      & \\
& \\
&\\
\texttt{Unit}   &No\\
&\\
\texttt{Comment}          & \\
&\\
\texttt{Ref.} &\cite{LodayRonco06} J.-L.~Loday, M. Ronco, \emph{On the structure of cofree Hopf}\\
& \emph{algebras}, J. reine angew. Math. 592 (2006), 123--155.\\
\end{tabular}

\newpage


\begin{tabular}{ll}
\texttt{Name}                   & {\bf Tridendriform  algebra}\index{algebra!tridendriform}\\
&\\
\texttt{Notation}               & $Tridend$\\
&\\
\texttt{Def.\ oper.}         & $x\l y,\  x\r y, \ x\cdot y$\\
\texttt{\hskip 1.2cm sym.}  & $ $\\
\texttt{\hskip 1.3cm rel.}    & 
$\left\{\begin{array}{rcl}
(x\l y)\l z &=& x\l (y*z)\  , \\
(x\r y)\l z &=& x\r (y\l z)\ , \\
(x * y)\r z &=& x\r (y\r z)\ ,\\
(x \r y) \cdot z &=& x \r (y \cdot z)\ ,\\
(x \l y) \cdot z &=& x \cdot (y \r z)\ ,\\
(x \cdot y) \l z &=& x \cdot (y \l z)\ , \\
(x \cdot y) \cdot z &=& x \cdot (y \cdot z)\ .
\end{array}
\right .
$\\
&where $x*y := x\l y+ x\r y + x\cdot y\ $. 
\\
&\\
\texttt{Free alg.}& Planar rooted trees with variables in between the leaves\\
&\\
\texttt{rep.~$\PP(n)$}      & $Tridend(n)=\KK[PT_n]\t \KK[\Sy_n],\  Tridend_n=\KK[PT_n]$\\
&\\
\texttt{$\dim \PP(n)$}      &$1, 3\times 2!, 11\times 3!, 45\times 4!, 197\times 5!,   903\times 6!, \ldots , C_{n}\times n! , \ldots$\\
 & where $C_n$ is the Schr\"oder (or super Catalan) number\\

&\\
\texttt{Gen.~series}& $f^{Tridend}(t)= \frac{-1+3t+\sqrt{1-6t+t^2}}{4t}$\\
&\\
\texttt{Dual operad}& $Triend^!=Trias$\\
&\\
\texttt{Chain-cplx}& Isomorphic to the total complex of a certain explicit tricomplex  \\
&\\
\texttt{Properties}          &  binary, quadratic, ns, Koszul.\\
&\\
\texttt{Alternative}         &  sometimes called dendriform trialgebras\\
& \\
\texttt{Relationship}      & $Tridend\alg \to As\alg,\quad  x*y:=x\l y + y \r x+ x\cdot y$,\\
&  $ComTridend\alg \to Tridend\alg,\quad  x\l y:=x\cdot y=: y \r x$\\
&\\
\texttt{Unit}       & $1\l x = 0 ,\ x\l 1 = x ,\quad
1\r x = x ,\  x\r 1 = 0, \quad 1\cdot x= 0 = x\cdot 1 $.\\
&\\
\texttt{Comment}          &There exist several variations (see \cite{Chapoton02} for instance).\\
&\\
\texttt{Ref.} &\cite{LodayRonco04} J.-L.\ Loday and M.\ Ronco, \emph{Trialgebras and families of polytopes,}\\
& Contemporary Mathematics (AMS) 346 (2004), 369--398.\\
\end{tabular}

\newpage


\begin{tabular}{ll}
\texttt{Name}                   & {\bf Triassociative algebra}\index{algebra!triassociative}\\
&\\
\texttt{Notation}               & $Trias$\\
&\\
\texttt{Def.\ oper.}         & $x\g y, \ x\d y,\  x\m y$\quad (no symmetry)\\
\texttt{\hskip 1.3cm rel.}    &
$\left\{\begin{array}{rcl}
(x\g y)\g z &=& x\g (y\g z), \\
(x\g y)\g z &=& x\g (y\d z), \\
(x\d y)\g z &=& x\d (y\g z) , \\
(x\g y)\d z &=& x\d (y\d z),\\
(x\d y)\d z &=& x\d (y\d z),\\
(x\g y)\g z &=& x\g (y\m z) , \\
(x\m y)\g z &=& x\m (y\g z) , \\
(x\g y)\m z &=& x\m (y\d z) , \\
(x\d y)\m z &=& x\d (y\m z) , \\
(x\m y)\d z &=& x\d (y\d z) , \\
(x\m y)\m z &=& x\m (y\m z) . 
\end{array}
\right .
$
\\
&one relation for each cell of the pentagon.\\
&\\
\texttt{Free alg.}&noncommutative polynomials with several variables marked \\
&\\
\texttt{rep.~$\PP(n)$}      & $Trias(n)=\KK^{2^n-1}\t \KK[\Sy_n], \quad Trias_n= \KK^{2^n-1}$\\
&\\
\texttt{$\dim \PP(n)$}    
&$1, 3\times 2!, 7\times 3!, 15\times 4!, 31\times 5!, 63\times 6!, 127\times 7!,\ldots , (2^n-1)\times n! , \ldots  $\\
&\\
\texttt{Gen.~series}& $ f^{Trias}(t)= \frac{t}{(1-t)(1-2t)} $\\
&\\
\texttt{Dual operad}& $Trias^!=Tridend$\\
&\\
\texttt{Properties}          &  binary, quadratic, ns, set-theoretic, Koszul.\\
&\\
\texttt{Alternative}         & Also called \emph{associative trialgebra}, or for short, \emph{trialgebra}.  \\
& \\
\texttt{Relationship}      &$As\alg \to Trias\alg,\quad  x\g y=x\d y=x\m y = xy$\\
&\\
\texttt{Unit}   &Bar-unit: $x\g 1=x=1\d x,\ 1\g x= 0 = x\d 1, \ 1\m x= 0 = x\m 1$\\
&\\
\texttt{Comment}          &Relations easy to understand in terms of planar trees \\
&\\
\texttt{Ref.}&\cite{LodayRonco04} J.-L.~Loday and M.~Ronco, \emph{Trialgebras and families of polytopes,}\\
&  Contemporary Mathematics (AMS) 346 (2004), 369--398.\\
\end{tabular}

\newpage


\begin{tabular}{ll}
\texttt{Name}                   & {\bf PostLie  algebra}\index{algebra!post-Lie}\\
&\\
\texttt{Notation}               & $x\circ y, \ [x, y]$\\
&\\
\texttt{Def.\ oper.}         & $x\circ y, \ [x, y]$\\
\texttt{\hskip 1.2cm sym.}  & $[x, y]=-[y, x]$\\
\texttt{\hskip 1.3cm rel.}    & $[x,[y,z]] + [y,[z,x]]+ [z,[x,y]] =0$\\
&$(x\circ y)\circ z - x\circ (y\circ z) - (x\circ z) \circ y + x\circ (z\circ y) = x \circ [y,z]$\\
&$[x,y]\circ z = [x\circ z, y] + [x,y\circ z]$\\
&\\
\texttt{Free alg.}&$PostLie(V)\cong Lie(??(V))$   \\
&\\
\texttt{rep.~$\PP(n)$}      & $ $\\
&\\
\texttt{$\dim \PP(n)$}       & $  1 , 3 , 20 , 210 ,3024, \cdots, a(n), \cdots $  \\
&\\
\texttt{Gen.~series}& $ f^{}(t)= -\log\big( \frac{1+\sqrt{1-4t}}{2})$\\
&\\
\texttt{Dual operad}& $PostLie^!=ComTrias$\\
&\\
\texttt{Chain-cplx}&   \\
&\\
\texttt{Properties}          &  binary, quadratic, Koszul.\\
&\\
\texttt{Alternative}         &  \\
& \\
\texttt{Relationship}      &$PostLie\alg \to ??\alg, \ xy = x\circ y ???$\\
&$PostLie\alg \to Lie\alg, \ \{x,y\} = x\circ y  - y\circ x +[x,y]$ \\
&$PreLie\alg\to PostLie\alg ,\ [x,y]=0$\\
&\\
\texttt{Unit}     &\\
&\\
\texttt{Comment}          & \\
&\\
\texttt{Ref.}&\cite{Vallette07} Vallette B., \emph{Homology of generalized partition posets}, \\
 & J. Pure Appl. Algebra 208 (2007), no. 2, 699--725.\\
\end{tabular}

\newpage


\begin{tabular}{ll}
\texttt{Name}                   & {\bf Commutative triassociative  algebra}\index{algebra!comutative triassociative}\\
&\\
\texttt{Notation}               & $ComTrias$\\
&\\
\texttt{Def.\ oper.}         & $x\g y, x\m y$\\
\texttt{\hskip 1.2cm sym.}  & $x\m y = y \m x$\\
\texttt{\hskip 1.3cm rel.}    & $\left\{\begin{array}{rcl}
(x\g y)\g z &=& x\g (y\g z), \\
(x\g y)\g z &=& x\g (z\g y), \\
(x\g y)\g z &=& x\g (y\m z) , \\
(x\m y)\g z &=& x\m (y\g z) , \\
(x\m y)\m z &=& x\m (y\m z) . 
\end{array}
\right .
$\\
&\\
\texttt{Free alg.}&  \\
&\\
\texttt{rep.~$\PP(n)$}      & $ $\\
&\\
\texttt{$\dim \PP(n)$}       & $ $\\
&\\
\texttt{Gen.~series}& $ f^{}(t)$\\
&\\
\texttt{Dual operad}& $ ComTrias^! = PostLie$\\
&\\
\texttt{Chain-cplx}&   \\
&\\
\texttt{Properties}          &  binary, quadratic, Koszul.\\
&\\
\texttt{Alternative}         & Triassociative with the following symmetry:\\
&$x\g y = y \d x$ and $x\m y = y \m x$\\
& \\
\texttt{Relationship}      & $ComTrias\alg \to Perm\alg$\\
&\\
\texttt{Unit}     &$x\g 1 = x, 1\g x = 0 , 1 \m x=0$ \\
&\\
\texttt{Comment}          & \\
&\\
\texttt{Ref.}&\\
\end{tabular}

\newpage


\begin{tabular}{ll}
\texttt{Name}                   & {\bf Commutative tridendriform  algebra}\index{algebra!commutative tridendriform}\index{algebra!CTD algebra}\\
&\\
\texttt{Notation}               & $CTD$\\
&\\
\texttt{Def.\ oper.}         & $x\l y,\ x\cdot y$\\
\texttt{\hskip 1.2cm sym.}  & $ x\cdot y=y\cdot x$\\
\texttt{\hskip 1.3cm rel.}    & 
$\left\{\begin{array}{rcl}
(x\l y)\l z &=& x\l (y\l z)+ x\l (z\l y)   , \\
(x \cdot y) \l z &=& x \cdot (y \l z)\ , \\
(x \l z) \cdot y &=& x \cdot (y \l z)\ , \\
(x \cdot y) \cdot z &=& x \cdot (y \cdot z)\ .
\end{array}
\right .
$
\\
&\\
\texttt{Free alg.}&$CTD(V)= $quasi-shuffle algebra on $V\ = QSym(V)$  \\
&\\
\texttt{rep.~$\PP(n)$}      & $CTD(n)= $\\
&\\
\texttt{$\dim \PP(n)$}      &$1, 3  , 13, 75, 541, 4683,  \ldots$\\
 &  \\

&\\
\texttt{Gen.~series}& $f^{CTD}(t)= \frac{\exp(t) -1}{2-\exp (t)} $\\
&\\
\texttt{Dual operad}& $CTD^!= $\\
&\\
\texttt{Chain-cplx}&   \\
&\\
\texttt{Properties}          &  binary, quadratic, Koszul.\\
&\\
\texttt{Alternative}         &  Handy to introduce $x*y:= x\l y + y \l x+ x\cdot y$ (assoc. and comm.) \\
& equivalently: tridendriform with symmetry:\\
&$x\l y = y \r x, \ x\cdot y = y \cdot x $\\
& \\
\texttt{Relationship}      &$CTD\alg \to Tridend\alg$\\
&\\
\texttt{Unit}       & $1\l x = 0 ,\ x\l 1 = x ,\quad 1\cdot x= 0 = x\cdot 1 $.\\
&\\
\texttt{Comment}          &\\
&\\
\texttt{Ref.} &\cite{Loday07} J.-L.~Loday, \emph{On the algebra of quasi-shuffles}, \\
& Manuscripta mathematica 123 (1), (2007), 79--93.\\
\end{tabular}

\newpage


\begin{tabular}{ll}
\texttt{Name}                   & {\bf Dual CTD  algebra}\index{algebra!dual commutative tridendriform}\index{algebra!dual CTD}\\
&\\
\texttt{Notation}               & $CTD^!$\\
&\\
\texttt{Def.\ oper.}         & $x\g y,\  [x, y]$ \\
\texttt{\hskip 1.2cm sym.}  & $[x, y]=-[y,x] $\\
\texttt{\hskip 1.3cm rel.}    & 
\\
&\\
\texttt{Free alg.}& \\
&\\
\texttt{rep.~$\PP(n)$}      &  \\
&\\
\texttt{$\dim \PP(n)$}      &$1, 3 , \ldots$\\
 &  \\

&\\
\texttt{Gen.~series}& $f^{CTD^!}(t)=  $\\
&\\
\texttt{Dual operad}& $(CTD^!)^!=CTD$\\
&\\
\texttt{Chain-cplx}& \\
&\\
\texttt{Properties}          &  binary, quadratic, Koszul.\\
&\\
\texttt{Alternative}         &   \\
&  \\
& \\
\texttt{Relationship}      &$Trias\alg\to CTD^!\alg  $\\
&\\
\texttt{Unit}       &  \\
&\\
\texttt{Comment}          & \\
&\\
\texttt{Ref.} & \\
\end{tabular}

\newpage


\begin{tabular}{ll}
\texttt{Name}                   & {\bf Gerstenhaber  algebra}\\
&\\
&underlying objects: graded vector spaces\\
&\\
\texttt{Notation}               & $Gerst$\\
&\\
\texttt{Def.\ oper.}         & $m=$ binary operation of degree $0$, $c=$ binary operation of degree $1$\\
\texttt{\hskip 1.2cm sym.}  & $m$ symmetric, $c$ antisymmetric\\
\texttt{\hskip 1.3cm rel.}    & $ c\circ_1 c + (c\circ_1 c)^{(123)}  + (c\circ_1 c)^{(321)}=0$, \\
&$c\circ_1 m - m\circ_2 c - (m\circ_1 c)^{(23)}=0$,\\
&$m\circ_1 m - m \circ_2 m=0$\ .\\
\texttt{Free alg.}&   \\
&\\
\texttt{rep.~$\PP(n)$}      & $ $\\
&\\
\texttt{$\dim \PP(n)$}       & $ $\\
&\\
\texttt{Gen.~series}& $ f^{ }(t)=$\\
&\\
\texttt{Dual operad}& $ $\\
&\\
\texttt{Chain-cplx}&   \\
&\\
\texttt{Properties}          &  binary, quadratic, Koszul.\\
&\\
\texttt{Alternative}         &   \\
& \\
\texttt{Relationship}      & \\
&\\
\texttt{Unit}     &\\
&\\
\texttt{Comment}          &To get the relations in terms of elements, do not forget to apply\\
& the Koszul sign rule. Observe that the last relation is associativity. \\
&\\
\texttt{Ref.}&\\
\end{tabular}

\newpage


\begin{tabular}{ll}
\texttt{Name}                   & {\bf Batalin-Vilkovisky  algebra}\\
&\\
&underlying objects: graded vector spaces\\
&\\
\texttt{Notation}               & $BV$\alg\\
&\\
\texttt{Def.\ oper.}         & $\DD$ unary degree $1$, $m$ binary degree $0$, $c$ binary degree $1$\\
\texttt{\hskip 1.2cm sym.}  & $m$ symmetric, $c$ antisymmetric\\
\texttt{\hskip 1.3cm rel.}    & $m\circ_1 m - m \circ_2 m=0$, \\
&$\Delta^2=0$, \\
&$c=\Delta\circ _1 m + m\circ_1 \Delta + m \circ_2  \Delta$,\\
&$c\circ_1 c + (c\circ_1 c)^{(123)}  + (c\circ_1 c)^{(321)}=0$, \\
&$c\circ_1 m - m\circ_2 c - (m\circ_1 c)^{(23)}=0$,\\
&$\Delta\circ _1 c + c \circ_1 \Delta +c \circ_2  \Delta=0$\ .\\
\texttt{Free alg.}&   \\
&\\
\texttt{rep.~$\PP(n)$}      & $ $\\
&\\
\texttt{$\dim \PP(n)$}       & $ $\\
&\\
\texttt{Gen.~series}& $ f^{ }(t)=$\\
&\\
\texttt{Dual operad}& $ $\\
&\\
\texttt{Chain-cplx}&   \\
&\\
\texttt{Properties}          &  unary and binary, inhomogeneous quadratic, Koszul.\\
&\\
\texttt{Alternative}         &Generated by $\DD$ and $m$ only.   \\
& \\
\texttt{Relationship}      &$BV\alg \to Gerst\alg$ \\
&\\
\texttt{Unit}     &\\
&\\
\texttt{Comment}          &Defined by J.-L. Koszul. To get the relations in terms of elements,\\
& do not forget to apply the Koszul sign rule. \\
\texttt{Ref.}&\cite{GVT09}\\
\end{tabular}

\newpage


\begin{tabular}{ll}
\texttt{Name}& {\bf   Magmatic algebra}\index{algebra!magmatic}\\
&\\
\texttt{Notation}& $Mag $\\
&\\
\texttt{Def.\ oper.}& $xy$\\
\texttt{\hskip 1.2cm sym.}  & $ $\\
\texttt{\hskip 1.3cm rel.} &$ $\\
&\\
\texttt{Free alg.}& any parenthesizing of words \\
&\\
\texttt{rep.~$\PP(n)$}& $Mag(n)=\KK[PBT_{n}]\t \KK[\Sy_n], \quad Mag_n= \KK[PBT_{n}]$\\
& (planar binary trees with $n$ leaves)\\
&\\
\texttt{$\dim \PP(n)$}
&$1, 1\times 2!, 2\times 3!, 5\times 4!, 14\times 5!, 42\times 6!, 132\times 7!, \ldots , c_{n-1}\times n!, \ldots $\\
&where $c_{n}=\frac{1}{n+1}\binom{2n}{n}$ (Catalan number)\\

&\\
\texttt{Gen.~series}& $f^{Mag}(t)= (1/2)(1-\sqrt{1-4t})$\\
&\\
\texttt{Dual operad}& $Mag^!=Nil_2$\\
&\\
\texttt{Chain-cplx}&   \\
&\\
\texttt{Properties}&  binary, quadratic, ns, set-theoretic, Koszul.\\
&\\
\texttt{Alternative}& sometimes (improperly) called \emph{nonassociative algebra}.   \\
& \\
\texttt{Relationship}& many ``inclusions" (all types of alg.~with only one gen.~op.)\\
&\\
\texttt{Unit}& $1x=x=x1$\\
&\\
\texttt{Comment}& \\
&\\
\texttt{Ref.}
\end{tabular}

\newpage


\begin{tabular}{ll}
\texttt{Name}                   & {\bf 2-Nilpotent  algebra}\index{algebra!nilpotent}\\
&\\
\texttt{Notation}               & $Nil_2$\\
&\\
\texttt{Def.\ oper.}         & $xy$\\
\texttt{\hskip 1.2cm sym.}  & $ $\\
\texttt{\hskip 1.3cm rel.}    & $(xy)z=0=x(yz) $\\
&\\
\texttt{Free alg.}&$Nil_2(V)= V\oplus V^{\t 2}$  \\
&\\
\texttt{rep.~$\PP(n)$}      & $Nil_2(2)= \KK[\Sy_2],\ Nil_2(n)= 0 $ for $n\geq 3$\\
&\\
\texttt{$\dim \PP(n)$}       & $1, 2, 0, 0, 0, 0, 0, \ldots  $\\
&\\
\texttt{Gen.~series}& $f^{\PP}(t)=t+t^2 $\\
&\\
\texttt{Dual operad}& ${Nil_2}^!=Mag$\\
&\\
\texttt{Chain-cplx}&   \\
&\\
\texttt{Properties}          &  binary, quadratic, ns, Koszul.\\
&\\
\texttt{Alternative}         &   \\
& \\
\texttt{Relationship}      & \\
&\\
\texttt{Unit}   &no\\
&\\
\texttt{Comment}          & \\
&\\
\texttt{Ref.}
\end{tabular}

\newpage


\begin{tabular}{ll}
\texttt{Name}& {\bf   Commutative magmatic algebra}\index{algebra!magmatic}\\
&\\
\texttt{Notation}& $ComMag $\\
&\\
\texttt{Def.\ oper.}& $x\cdot y$\\
\texttt{\hskip 1.2cm sym.}  & $x\cdot y=y\cdot x $\\
\texttt{\hskip 1.3cm rel.} &none\\
&\\
\texttt{Free alg.}& any parenthesizing of commutative words \\
&\\
\texttt{rep.~$\PP(n)$}& $ComMag(n)=\KK[shBT(n)]$\\
&$shBT(n)=$ \{shuffle binary trees with $n$ leaves\}\\
&\\
\texttt{$\dim \PP(n)$}&$\dim ComMag(n)= (2n-3)!! = 1\times 3\times \cdots \times (2n-3)$\\
&\\
\texttt{Gen.~series}& $f^{ComMag}(t)= $\\
&\\
\texttt{Dual operad}& $ComMag^!\alg$: $[x,y]$ antisymmetric, $[[x,y],z]=0$.\\
&\\
\texttt{Chain-cplx}&   \\
&\\
\texttt{Properties}&  binary, quadratic, set-theoretic, Koszul.\\
&\\
\texttt{Alternative}&   \\
& \\
\texttt{Relationship}& $ComMag\mono PreLie, x\cdot y := \{x,y\} + \{y,x\}$, cf.\ \cite{BergeronLoday11}\\
&\\
\texttt{Unit}& $1x=x=x1$\\
&\\
\texttt{Comment}& \\
&\\
\texttt{Ref.}
\end{tabular}

\newpage


\begin{tabular}{ll}
\texttt{Name}& {\bf   Anti-symmetric nilpotent algebra}\index{algebra!anti-symmetric nilpotent}\\
&\\
\texttt{Notation}& $ComMag^{!} $\\
&\\
\texttt{Def.\ oper.}& $x\cdot y$\\
\texttt{\hskip 1.2cm sym.}  & $x\cdot y=-y\cdot x $\\
\texttt{\hskip 1.3cm rel.} &$(xy)z=0$\\
&\\
\texttt{Free alg.}&   \\
&\\
\texttt{rep.~$\PP(n)$}&  \\
& \\
&\\
\texttt{$\dim \PP(n)$}& \\
&\\
\texttt{Gen.~series}&  \\
&\\
\texttt{Dual operad}& $ComMag $.\\
&\\
\texttt{Chain-cplx}&   \\
&\\
\texttt{Properties}&  binary, quadratic, set-theoretic, Koszul.\\
&\\
\texttt{Alternative}&   \\
& \\
\texttt{Relationship}&  \\
&\\
\texttt{Unit}&  \\
&\\
\texttt{Comment}& \\
&\\
\texttt{Ref.}
\end{tabular}

\newpage

\begin{tabular}{ll}
\texttt{Name}                   & {\bf  Quadri-algebra}\index{algebra!quadri-}\\
&\\
\texttt{Notation}               & $Quadri$\\
&\\
\texttt{Def.\ oper.}         & $x\nw y,\ x\ne y,\ x\se y,\ x\sw y$\quad  called NW, NE, SE, SW oper.\\
\texttt{\hskip 1.3cm rel.}    & \\
\end{tabular}
\begin{equation*}  \label{E:quadri}
\raisebox{25pt}{\xymatrix@C=5pt@R=5pt{
{\mbox{\footnotesize $(x\nw y)\nw z =x\nw(y\star z)$}} & &
{\mbox{\footnotesize $(x\ne y)\nw z = x\ne(y\west z)$}} &  &
{\mbox{\footnotesize $(x\north y)\ne z=x\ne(y\east z)$}}  \\
 {\mbox{\footnotesize $(x\sw y)\nw z =x\sw(y\north z)$}} & &
 {\mbox{\footnotesize $(x\se y)\nw z = x\se(y\nw z)$}} &  &
{\mbox{\footnotesize $(x\south y)\ne z=x\se(y\ne z)$}}             \\
 {\mbox{\footnotesize $(x\west y)\sw z=x\sw(y\south z)$}} & &
 {\mbox{\footnotesize $(x\east y)\sw z = x\se(y\sw z)$}} &  &
{\mbox{\footnotesize $(x\star y)\se z =x\se(y\se z)$}}
}}
\end{equation*}
\begin{tabular}{ll} where&\\
&$x\east y:= x\ne y+x\se y, \quad x\west y:= x\nw y+x\sw y$\\
&$x\south y:= x\se y+x\sw y, \quad x\north y:= x\ne y+x\nw y$\\
&$x\star y:=x\se y+x\ne y+x\nw y+x\sw y $\\
&$\qquad \quad = x\east y+x\west y=x\south y+x\north y$\\
&\\
\texttt{Free alg.}&   \\
&\\
\texttt{rep.~$\PP(n)$}      & $ $\\
&\\
\texttt{$\dim \PP(n)$}       & $ 1, 4\times 2!, 23\times 3!, 156\times 4!, 1162\times 5!, 9192\times 6!, \ldots$\\
&$\dim \PP_n =  \frac{1}{n} \sum_{j=n}^{2n-1}{3n \choose n+1+j}{j-1 \choose j-n}$\\
\texttt{Gen.~series}& $ f^{}(t)$\\
&\\
\texttt{Dual operad}& $Quadri^!$\\
&\\
\texttt{Properties}          &  binary, quadratic, ns, Koszul.\\
&\\
\texttt{Alternative}         & $Quadri = Dend\ \blacksquare\  Dend= PreLie\bullet Dend$   \\
& \\
\texttt{Relationship}      & Related to dendriform in several ways\\
&\\
\texttt{Unit}     &partial unit (like in dendriform)\\
&\\
\texttt{Comment}          &There exist several variations like $\PP \blacksquare \PP$ (cf.\ Ph.\ Leroux \cite{Leroux04})  \\
&\\
\texttt{Ref.}&\cite{AguiarLoday04} M.\ Aguiar, J.-L.\ Loday, \emph{Quadri-algebras},\\
& J. Pure Applied Algebra 191 (2004), 205--221.\\
\end{tabular}


\begin{tabular}{ll}
\texttt{Name}                   & {\bf Dual quadri-algebra}\index{algebra!dual quadri-}\\
&\\
\texttt{Notation}               & $Quad^!$\\
&\\
\texttt{Def.\ oper.}         &  $x\nw y,\ x\ne y,\ x\se y,\ x\sw y$\\
\texttt{\hskip 1.2cm sym.}  & $ $\\
\texttt{\hskip 1.3cm rel.}    & A FAIRE  $ $\\
&\\
\texttt{Free alg.}&   \\
&\\
\texttt{rep.~$\PP(n)$}      & \\
& \\
\texttt{$\dim \PP(n)$}       &  $ 1, 4\times 2!, 9\times 3!, 16\times 4!, 25\times 5!, \ldots , n^2\times n!, \ldots $\\
&\\
\texttt{Gen.~series}& $ f^{}(t)$\\
&\\
\texttt{Dual operad}& $ $\\
&\\
\texttt{Chain-cplx}&   \\
&\\
\texttt{Properties}          &  binary, quadratic, Koszul.\\
&\\
\texttt{Alternative}         &   \\
& \\
\texttt{Relationship}      & \\
&\\
\texttt{Unit}     &\\
&\\
\texttt{Comment}          & \\
&\\
\texttt{Ref.}&\cite{Vallette08}  B.Vallette, \emph{Manin products, Koszul duality,}\\
& \emph{Loday algebras and Deligne conjecture}.\\
& J. Reine Angew. Math. 620 (2008), 105--164.\\
\end{tabular}


\begin{tabular}{ll}
\texttt{Name}                   & {\bf Duplicial  algebra}\index{algebra!duplicial}\index{algebra!duplicial}\\
&\\
\texttt{Notation}               & $Dup$\\
&\\
\texttt{Def.\ oper.}         & $x\l y,\  x\r y$\\
\texttt{\hskip 1.2cm sym.}  & $ $\\
\texttt{\hskip 1.3cm rel.}    & 
$\left\{\begin{array}{rcl}
(x\l y)\l z &=& x\l (y\l z)   , \\
(x\r y)\l z &=& x\r (y\l z) , \\
 (x\r y)\r z&=& x\r (y\r z).
\end{array}
\right .
$
\\
&\\
\texttt{Free alg.}& $Dup(V)=\oplus_{n\geq 1} \KK[PBT_{n+1}]\t V^{\t n}$,  for p.b.~trees $s$ and $t$: \\
& $x\r y$ \emph{over} operation is grafting of $x$ on the leftmost leaf of $y$\\
& $x\l y$ \emph{under} operation is grafting of $y$ on the rightmost leaf of $x$\\
&\\
\texttt{rep.~$\PP(n)$}      & $Dup(n)=  \KK[PBT_{n+1}]\t \KK[\Sy_n] $, so $Dup_n=  \KK[PBT_{n+1}]$\\
&\\
\texttt{$\dim \PP(n)$}   &$1, 2\times 2!, 5\times 3!, 14\times 4!, 42\times 5!, 132\times 6!, 429\times 7!, \ldots , c_{n}\times n!, \ldots$\\
 &where $c_n=\frac{1}{n+1}\binom{2n}{n}$ is the Catalan number\\

&\\
\texttt{Gen.~series}& $f^{Dup}(t)= \frac{1-2t-\sqrt{1-4t}}{2t}=y, \quad y^2-(1-2t)y+t=0$\\
&\\
\texttt{Dual operad}& $Dup^!$\\
&\\
\texttt{Chain-cplx}& Isomorphic to the total complex of a certain explicit bicomplex  \\
&\\
\texttt{Properties}          &  binary, quadratic, ns, set-theoretic, Koszul.\\
&\\
\texttt{Alternative}         &   \\
& \\
\texttt{Relationship} & $2as\alg \to Dup\alg \to As^2\alg$\\
&\\
\texttt{Unit}       & \\
&\\
\texttt{Comment}          &The associator of $xy:=x\r y - x\l y $ is \\
& $as(x,y,z) = x\l (y\r z)-(x\l y)\r z$. \\
&Appeared first in \cite{BrouderFrabetti03}\\
\texttt{Ref.} &\cite{Loday08} J.-L. Loday, \emph{Generalized bialgebras and triples of operads},\\
& Ast\' erisque (2008), no 320, x+116 p.\\
\end{tabular}

\newpage


\begin{tabular}{ll}
\texttt{Name}                   & {\bf Dual duplicial  algebra}\index{algebra!dual duplicial}\index{algebra!dual duplicial}\\
&\\
\texttt{Notation}               & $Dup^!$\\
&\\
\texttt{Def.\ oper.}         & $x\l y, \ x\r y$\\
\texttt{\hskip 1.2cm sym.}  & $ $\\
\texttt{\hskip 1.3cm rel.}    & $\left\{\begin{array}{rcl}
(x\l y)\l z &=& x\l (y\l z), \\
(x\l y)\r z &=&0, \\
(x\r y)\l z &=& x\r (y\l z) , \\
               0&=& x\l (y\r z),\\
(x\r y)\r z &=& x\r (y\r z).
\end{array}
\right .
$
\\
&\\
\texttt{Free alg.}&$Dup^!(V)=\oplus_{n\geq 1}\ n\ \Vtn$ \\
&noncommutative polynomials with one variable marked  \\
&\\
\texttt{rep.~$\PP(n)$}      & $Dup^!(n)=\KK^n\t \KK[\Sy_n], \quad Dup^!_n= \KK^n$\\
&\\
\texttt{$\dim \PP(n)$}       &$1, 2\times 2!, 3\times 3!, 4\times 4!, 5\times 5!, 6\times 6!, 7\times 7!, \ldots , n\times n!, \ldots,$\\
&\\
\texttt{Gen.~series}& $ f^{Dup^!}(t)= \frac{t}{(1-t)^2}$\\
&\\
\texttt{Dual operad}& $Dup^{!!}=Dup$\\
&\\
\texttt{Chain-cplx}& see ref.  \\
&\\
\texttt{Properties}          &  binary, quadratic, Koszul.\\
&\\
\texttt{Alternative}         &   \\
& \\
\texttt{Relationship}     
&\\
\texttt{Unit}   & \\
&\\
\texttt{Comment}          & \\
&\\
\texttt{Ref.}& \\
\end{tabular}

\newpage


\begin{tabular}{ll}
\texttt{Name}                   & {\bf  $As^{(2)}$-algebra}\\
&\\
\texttt{Notation}               & $As^{(2)}$\\
&\\
\texttt{Def.\ oper.}         & $x*y, x\cdot y$\\
\texttt{\hskip 1.2cm sym.}  & $ $\\
\texttt{\hskip 1.3cm rel.}    & $(x\circ_{1}y)\circ_{2}z = x\circ_{1}(y\circ_{2}z)$ any values (4 relations)\\
&\\
\texttt{Free alg.}&   \\
&\\
\texttt{rep.~$\PP(n)$}      & $As^{(2)}_{n}= \KK[\{0,1\}^{n-1}]$\\
&\\
\texttt{$\dim \PP(n)$}       & $2^{n-1}$\\
&\\
\texttt{Gen.~series}& $ f^{As^{(2)}}(t)=\frac{t}{1-2t}$\\
&\\
\texttt{Dual operad}& ${As^{(2)}}^!=As^{(2)}$\\
&\\
\texttt{Chain-cplx}&   \\
&\\
\texttt{Properties}          &  binary, quadratic, set-theoretic, ns.\\
&\\
\texttt{Alternative}         &   \\
& \\
\texttt{Relationship}      &$2as\alg \to Dup\alg \to As^{(2)}\alg$ \\
&\\
\texttt{Unit}     &\\
&\\
\texttt{Comment}          & Variations: $ As^{(k)}$\\
&\\
\texttt{Ref.}&\\
\end{tabular}

\newpage


\begin{tabular}{ll}
\texttt{Name}                   & {\bf  $As^{\langle2\rangle}$-algebra}\\
&\\
\texttt{Notation}               & $As^{\langle2\rangle}$, compatible products algebra\\
&\\
\texttt{Def.\ oper.}         & $x*y, x\cdot y$\\
\texttt{\hskip 1.2cm sym.}  & $ $\\
\texttt{\hskip 1.3cm rel.}    & $(x*y)*z = x*(y*z)$\\
  & $(x*y)\cdot z + (x\cdot y)*z  = x*(y\cdot z) + x\cdot (y*z)$\\
   & $(x\cdot y)\cdot z = x\cdot (y\cdot z)$\\
&\\
\texttt{Free alg.}& similar to dendriform  \\
&\\
\texttt{rep.~$\PP_{n}$}      & $As^{\langle2\rangle}_{n}= \KK[PBT_{n+1}]$\\
&\\
\texttt{$\dim \PP_{n}$}       & $c_n=\frac{1}{n+1}\binom{2n}{n}$\\
&\\
\texttt{Gen.~series}& $ f^{As^{\langle2\rangle}}(t)= \frac{1-2t-\sqrt{1-4t}}{2t}$\\
&\\
\texttt{Dual operad}&  An $(As^{\langle 2\rangle})^!$-algebra has 2 generating binary operations and 5 relations:\\
& $(x\circ_{i}y)\circ_{j}z = x\circ_{i}(y\circ_{j}z) $ and\\
& $(x*y)\cdot z + (x\cdot y)*z  = x*(y\cdot z) + x\cdot (y*z)$\\
&\\
\texttt{Chain-cplx}&   \\
&\\
\texttt{Properties}          &  binary, quadratic, ns, Koszul.\\
&\\
\texttt{Alternative}         &   \\
& \\
\texttt{Relationship}      &$(As^{\langle 2\rangle})^!\alg \to As^{(2)}\alg \to As^{\langle 2\rangle}\alg$ \\
&\\
\texttt{Unit}     &\\
&\\
\texttt{Comment}          &equivalently $\lambda\ x*y + \mu\ x\cdot y$ is associative for any $\lambda, \mu$\\
& Variations: $As^{\langle k \rangle}$, $Hoch$-alg.\\
&\\
\texttt{Ref.}& for instance \cite{Goncharov05}, \cite{OdesskiiSokolov06}, \cite{Dotsenko09}.\\
\end{tabular}

\newpage


\begin{tabular}{ll}
\texttt{Name}                   & {\bf  L-dendriform algebra}\\
&\\
\texttt{Notation}               & $L\textrm{-}dend$\\
&\\
\texttt{Def.\ oper.}         &  $ x\vartriangleright y $ and $ x\vartriangleleft y$ \\
\texttt{\hskip 1.2cm sym.}  &\\
\texttt{\hskip 1.3cm rel.}    & $ x\vartriangleright(y \vartriangleright z) -(x \bullet y) \vartriangleright z = y \vartriangleright(x \vartriangleright z)-(y \bullet x) \vartriangleright z$\\
&$x\vartriangleright(y \vartriangleleft z) - (x \vartriangleright y) \vartriangleleft z = y \vartriangleleft (x \bullet z) - (y \vartriangleleft x)\vartriangleleft z$\\
& where $x \bullet y := x \vartriangleright y + x \vartriangleleft  y.$\\
\texttt{Free alg.}&   \\
&\\
\texttt{rep.~$\PP(n)$}      & $ $\\
&\\
\texttt{$\dim \PP(n)$}       & $ $\\
&\\
\texttt{Gen.~series}& $ f^{}(t)=$\\
&\\
\texttt{Dual operad}& $ $\\
&\\
\texttt{Chain-cplx}&   \\
&\\
\texttt{Properties}          & $L\textrm{-}dend = PreLie \bullet PreLie$ ? \\
&\\
\texttt{Alternative}         &   \\
& \\
\texttt{Relationship}      &$Dend\alg\to L\textrm{-}dend\alg$ via $x\vartriangleright y:= x\r y$ and $x\vartriangleleft y:= x\l y$ \\
&$L\textrm{-}dend\alg\to PreLie\alg, (A,\vartriangleright, \vartriangleleft)\mapsto (A, \bullet) $\\
\texttt{Unit}     &\\
&\\
\texttt{Comment}          &Various variations like $L\textrm{-}quad\alg$ \cite{BGN10} \\
&\\
\texttt{Ref.}& Bai C.M., Liu L.G., Ni X., \emph{L-dendriform algebras}, preprint (2010). \\
\end{tabular}

\newpage


\begin{tabular}{ll}
\texttt{Name}& {\bf   Lie-admissible algebra}\index{algebra!Lie-admissible}\\
&\\
\texttt{Notation}& $Lie\-adm$\\
&\\
\texttt{Def.\ oper.}& $xy$\\
\texttt{\hskip 1.2cm sym.}  & $ $\\
\texttt{\hskip 1.3cm rel.} &$[x,y]=xy-yx$ is a Lie bracket,  that is \\
&$\sum_{\ss} \sgn (\ss) \ss\big( (xy)z - x(yz)\big)=0$\\
&\\
\texttt{Free alg.}&   \\
&\\
\texttt{rep.~$\PP(n)$}& $Lie\-adm(n)= ?$\\
&\\
\texttt{$\dim \PP(n)$}
&$1, 2, 11, \ldots $\\
&\\
\texttt{Gen.~series}& $f^{Lie\-adm}(t)= ?$  its dimensions are $1,2,1, ??$\\
&\\
\texttt{Dual operad}& $Lie\-adm^!$\\
&\\
\texttt{Chain-cplx}&   \\
&\\
\texttt{Properties}&  binary, quadratic, Koszul  ??.\\
&\\
\texttt{Alternative}&   \\
& \\
\texttt{Relationship}&$As\alg \to PreLie\alg \to Lie\-adm\alg$\\
&\\
\texttt{Unit}&  \\
&\\
\texttt{Comment}& \\
&\\
\texttt{Ref.}
\end{tabular}

\newpage


\begin{tabular}{ll}
\texttt{Name}                   & {\bf PreLiePerm  algebra}\index{algebra!pre-Lie permutative}\index{algebra!PreLiePerm}\\
&\\
\texttt{Notation}               & $PreLiePerm$\\
&\\
\texttt{Def.\ oper.}         & $x\l y,\  x\r y,\  x*y$\\
\texttt{\hskip 1.2cm sym.}  & $$\\
\texttt{\hskip 1.3cm rel.}    & $\left\{\begin{array}{rcl}
(x\l y)\l z &=& x\l (y* z), \\
(x\l y)\l z &=& x\l (z* y), \\
(x* y)\r z &=& x\r (y\r z) , \\
(x\r y)\l z &=& x\r (z\r y) , \\
(x\r y)\l z &=& x\r (y\l z), 
\end{array}
\right .
$\\
\texttt{Free alg.}& $$\\
&\\
\texttt{rep.~$\PP(n)$}      & $$\\
&\\
\texttt{$\dim \PP(n)$}   
&$1, 6, $\\
&\\
\texttt{Gen.~series}& $$\\
&\\
\texttt{Dual operad}& $PermPreLie= Perm \circ PreLie=Perm\t PreLie$\\
&\\
\texttt{Chain-cplx}&   \\
&\\
\texttt{Properties}          &  binary, quadratic, set-theoretic.\\
&\\
\texttt{Alternative}         &   \\
& \\
\texttt{Relationship}      &$PreLiePerm= PreLie\bullet Perm$, see \cite{Vallette08}\\
&$Zinb\alg\to PreLiePerm\alg \to Dend\alg$ \\
&\\
\texttt{Unit}     &no\\
&\\
\texttt{Comment}          & \\
&\\
\texttt{Ref.}&\cite{Vallette08}  B.Vallette, \emph{Manin products, Koszul duality,}\\
& \emph{Loday algebras and Deligne conjecture}.\\
& J. Reine Angew. Math. 620 (2008), 105--164.\\
\end{tabular}%

\newpage

Integer sequences which appear in this paper, up to some shift and  up to multiplication by $n!$ or $(n-1)!$.

\bigskip
%
\\

\newpage


\bibliographystyle{amsalpha}
\bibliography{Encyc-bib}

\end{document}